\documentclass{article}
\usepackage[a4paper, margin=3cm]{geometry}
\usepackage[utf8]{inputenc}
\usepackage[english]{babel}

\usepackage{soul}
\usepackage{amsthm}
\usepackage{amssymb}
\usepackage{indentfirst}
\usepackage{amsmath}
\usepackage{graphicx, color}
\usepackage{enumitem}
\usepackage{algorithm}
\usepackage{algpseudocode}
\usepackage{authblk}

\usepackage{tikz}

\graphicspath{{Imagenes/}}

\newtheorem{theorem}{Theorem}[section]
\newtheorem{lemma}[theorem]{Lemma}
\newtheorem{corollary}[theorem]{Corollary}
\newtheorem{proposition}[theorem]{Proposition}

\newcommand{\exact}[2]{#1^{[\sharp #2]}}

\newcommand{\CV}{\mathrm{CV}}

\floatname{algorithm}{Procedure}

\title{Characterizing and recognizing exact-distance squares of graphs}

\author[1,2]{Yandong Bai} 
\author[3]{Pedro P. Cortés} 
\author[4]{Reza Naserasr}
\author[5]{Daniel A. Quiroz}

\affil[1]{\small Research $\&$ Development Institute of Northwestern Polytechnical University in Shenzhen, Shenzhen, Guangdong 518057, China}
\affil[2]{\small School of Mathematics and Statistics, Northwestern Polytechnical University, Xi'an, Shaanxi 710129, China, Email: \texttt{bai@nwpu.edu.cn}}
\affil[3]{\small Departamento de Ingenier\'ia Matem\'atica, Universidad de Chile, Santiago, Chile, Email: \texttt{pedrocortes1997@gmail.com}}
\affil[4]{\small Université Paris Cité, CNRS, IRIF, F-75013, Paris, France	
	Email: \texttt{reza@irif.fr}}	
\affil[5]{\small Instituto de Ingenier\'ia Matem\'atica-CIMFAV, Universidad de Valpara\'iso, Valpara\'iso, Chile. Email: \texttt{daniel.quiroz@uv.cl}}

\date{}

\begin{document}
	
	\maketitle
	
	\begin{abstract}
		For a graph $G=(V,E)$, its \emph{exact-distance square}, $\exact{G}{2}$, is the graph with vertex set~$V$ and with an edge between vertices $x$ and $y$ if and only if $x$ and $y$ have distance (exactly) $2$ in~$G$. The graph $G$ is an \emph{exact-distance square root} of   $\exact{G}{2}$. We give a characterization of graphs having an exact-distance square root, our characterization easily leading to a polynomial-time recognition algorithm. We show that it is NP-complete to recognize graphs with a bipartite exact-distance square root. These two results strongly contrast known results on (usual) graph squares. We then characterize graphs having a tree as an exact-distance square root, and from this obtain a polynomial-time recognition algorithm for these graphs. Finally, we show that, unlike for usual square roots, a graph might have (arbitrarily many) non-isomorphic exact-distance square roots which are trees. 
		
	\end{abstract}

	\section{Introduction}
	The notion of graph power is a fundamental one in graph theory. Given a graph $G$ and a positive integer $p$, the \emph{power graph} $G^p$ is the graph with vertex set $V$ and with an edge between vertices $x$ and $y$ if and only if $x$ and $y$ have distance at most $p$ in $G$. We say that $G^2$ is the \emph{square} of $G$, and that $G$ is a \emph{square root} of $G^2$. Since the introduction of the concept of powers, a key issue has been to characterize and efficiently recognize graphs which have square root. 	
	In 1967, Mukhopadhyay \cite{mukhopadhyay1967square} characterized graphs having at least one square root. This result was generalized by Geller \cite{geller1968square} who characterized digraphs having at least one square root. Yet this characterization does not lead to an efficient recognition algorithm. Moreover, in 1994, Motwani and Sudan \cite{motwani1994computing} proved that the problem of deciding whether a graph has a square root is NP-complete.
	
	In the paper in which the notion of square was first introduced, Ross and Harary \cite{ross1960square} characterized squares of trees. This characterization does lead to a polynomial-time (in fact, linear-time) algorithm to recognize squares of trees, as proved by Lin and Skiena \cite{lin1995algorithms} (see Lau \cite{lau2006bipartite} for a different such algorithm). Going beyond trees, for any graph class $\mathcal{C}$, the \emph{$\mathcal{C}$-root-problem} has as instance a graph $G$, and as task that of deciding whether $G$ has a square root $H$ which satisfies $H\in\mathcal{C}$. This problem has been amply studied \cite{ducoffe2019finding, farzad2012square, farzad2012complexity, golovach2019algorithms, lau2004recognizing}. Particularly interesting is the result of Lau~\cite{lau2006bipartite} that the problem is in P when $\mathcal{C}$ is the class of bipartite graphs; it had been conjectured by Motwani and Sudan \cite{motwani1994computing} that this problem would be NP-complete. Also remarkable is the recent result of Dvo\v{r}\'ak, Lahiri, and Moore \cite{dvovrak2022square} that the problem is NP-complete when $\mathcal{C}$ is the class of 6-apex graphs, and thus it can be NP-complete on sparse graph classes.
	
	In this paper we consider the following refinement of the notion of graph power. For a graph $G=(V,E)$ and a positive integer $p$, the \emph{exact distance-$p$ graph} $\exact{G}{p}$ is the graph with vertex set~$V$ and with an edge between vertices $x$ and $y$ if and only if $x$ and $y$ have distance (exactly) $p$ in~$G$. In other words, $\exact{G}{p}$ is obtained from $G^p$ by removing the edges of $G^{p-1}$. This notion has received increased attention in recent years, mostly from a graph coloring perspective \cite{foucaud2021exact, bousquet2019exact, quiroz2020colouring, la20222, almulhim2020estimating}, but also through algorithmic \cite{panda2022exact} and structural \cite{brevsar2019exact, foucaud2021cliques} perspectives. A result \cite{van2019chromatic, nevsetvril2012sparsity} which has sparked much of this research is the following. Let $\mathcal{C}$ be a (proper) minor closed class of graphs. For every integer $p\ge 1$, there exists a constant $N=N(\mathcal{C},p)$ such that $\chi(\exact{G}{p})\le N$ if $p$ is odd and  $\chi(\exact{G}{p})\le N \cdot \Delta(G)$ if $p$ is even. This result actually extends to any class with bounded expansion, and gives a more refined view on known results on the chromatic numbers of graph powers.
	
	We say that $\exact{G}{2}$ is the \emph{exact-distance square} of $G$, and that $G$ is an \emph{exact-distance square root} of   $\exact{G}{2}$.
	We study the characterization of graphs which have an exact-distance square root within a particular graph class. Moreover, we study the \emph{$\mathcal{C}$-exact-distance-root} problem, which has as instance a graph $G$, and as task that of deciding if $G$ has an exact-distance square root $H$ which satisfies $H\in\mathcal{C}$.
	
	We first consider the class of all graphs, and obtain the following characterization. Here we let $\overline{G}$ denote the complement of a graph $G$.
	
	\begin{theorem}\label{theorem:all}
		A graph $G$ has an exact-distance square root if and only if $G$ is equal, as a labelled graph, to $\overline{G}^{[\sharp 2]}$.
	\end{theorem}
	
	This characterization immediately gives a quadratic-time algorithm for the $\mathcal{C}$-exact-distance-root problem, when $\mathcal{C}$ is the class of all graphs (simply, compute the exact-distance square of the complement of $G$ and check the equality). Moreover, this algorithm can additionally output an exact-distance square root if there is one, at no additional computational cost.
	
	Together with the previous result, the following helps to grasp the strong contrast between squares and exact-distance squares.
	
	\begin{theorem}\label{theorem:bipartite}
		The $\mathcal{C}$-exact-distance-root problem is NP-complete when $\mathcal{C}$ is the class of bipartite graphs. 
	\end{theorem}
	
	We also characterize graphs which are the exact-distance square of some bipartite graph (see Theorem~\ref{thm:Clique-dualCharacterization}), and graphs which are the exact-distance square of some triangle-free graph (see Theorem~\ref{thm:Triangle-free-Root}). Whether the $\mathcal{C}$-exact-distance-root problem is NP-complete or not for the class of triangle-free graphs, we leave as an open problem.
	
	We then turn our attention to exact-distance squares of trees. We characterize graphs which are exact-distance squares of trees (see Corollary~\ref{corollary:trees}) and, based on this characterization, provide the following result.
	
	\begin{theorem}\label{theorem:trees}
		There is a polynomial-time algorithm for the $\mathcal{C}$-exact-distance-root problem when $\mathcal{C}$ is the class of all trees. 
	\end{theorem}
	
	In fact, this algorithm can be extended to obtain a tree exact-distance square root if there is one, also in polynomial time (see Section \ref{sec:recognition}).
	
	
	We end the paper with another result that contrasts with known result on squares. Ross and Harary \cite{ross1960square}, and later Lau \cite{lau2006bipartite},  proved that tree square roots, when they exists, are unique up to isomorphism. Farzad, Lau, Le and Tuy \cite{Farzadetal} showed that this holds beyond tree square roots to square roots with girth at least 7, and Adamaszek and Adamaszek \cite{adamasquared} further showed that this holds for square roots with girth at least 6, a result which is best possible. We show that for exact-distance squares this does not hold even for trees, that is, a graph might have (arbitrarily many) non-isomorphic exact-distance square roots which are trees.
	
	\begin{theorem}\label{theorem:unique}
		For every $m\ge 2$, there is a graph with at least $m!$ nonisomorphic exact-distance square roots which are trees.
	\end{theorem}
	
	The rest of the paper is organized as follows. In Section~\ref{sec:all} we prove Theorem~\ref{theorem:all} and show how to extend it to digraphs. In Section~\ref{sec:bipartite} we prove Theorem~\ref{theorem:bipartite}, and give a characterization of exact-distance squares of bipartite graphs and a characterization of exact-distance squares of triangle-free graphs. In Section~\ref{sec:charactrees} we characterize exact-distance squares of trees, and in Section~\ref{sec:recognition} we use this characterization to prove Theorem~\ref{theorem:trees}. Finally, in Section~\ref{sec:unique} we prove Theorem~\ref{theorem:unique}.
	
	We now clarify some notation and definitions. We consider only finite and simple graphs.  A \emph{block} in a graph $G$ is a maximal subgraph without cut-vertices. A connected graph is a said to be \emph{clique-tree} if each of its blocks is a complete graph,  equivalently, if each cycle induces a complete graph.

	\section{Exact-distance squares of (di)graphs}\label{sec:all}

	In this section we prove Theorem \ref{theorem:all} which characterizes graphs which have some exact-distance square root. After this, we mention how this result extends to digraphs.

	\begin{proof}[Proof of Theorem~\ref{theorem:all}]
		The sufficiency is obvious. It remains to show the necessity. Assume that $G$ has an exact-distance square root $H$, i.e., $G=\exact{H}{2}$. Note that $V(G)=V(H)=V(\overline{G})$.  Since an adjacent pair of vertices in $H$ are at distance 1, and not 2, they are not adjacent in $G$, thus $H\subseteq \overline{G}$. Similarly, since an adjacent pair of vertices in $\exact{\overline{G}}{2}$ are not adjacent in $\overline{G}$, we have $\exact{\overline{G}}{2} \subseteq G$. For a pair $x,y$ of adjacent vertices in $G$, since $G=\exact{H}{2}$, they must be at distance exactly 2 in $H$ and, by the definition of the complement, they must be at distance at least 2 in $\exact{\overline{G}}{2}$. But since $H\subseteq \overline{G}$, they should be at distance exactly 2 in $\exact{\overline{G}}{2}$ as well. Thus $ G \subseteq \exact{\overline{G}}{2}$. This implies $G=\exact{\overline{G}}{2}$.
	\end{proof}
	
	Following the same proof, Theorem~\ref{theorem:all} can be extended to digraphs with the following (natural) extensions of the definitions. In a digraph $D$ the \emph{distance} of the ordered pair $(x,y)$ is the length of a shortest directed path starting at $x$ and ending at $y$. When there is no such a path, the distance of $(x,y)$ is $\infty$. The \emph{exact-distance square} of a digraph $D$ is a digraph on the same set of vertices, with $(x,y)$ being an arc if and only if the distance of $(x,y)$ is 2 in $D$.  The \emph{complement} of a digraph $D$ contains all arcs $(x,y)$, $x\neq y$, which are not arcs of $D$. 
	
	\section{Exact-distance squares of bipartite graphs}\label{sec:bipartite}
	
	The result from the previous section tells us that for the class $\mathcal{C}$ of all graphs, while the $\mathcal{C}$-root problem is NP-hard, the $\mathcal{C}$-exact-distance-root problem can be solved in quadratic time. Here we show that the roles are reversed when we consider the class $\mathcal{B}$ of all bipartite graphs. For the $\mathcal{B}$-root problem,  Lau~\cite{lau2006bipartite} provided a polynomial-time algorithm. In contrast, we show that $\mathcal{B}$-exact-distance-root problem is among NP-complete problems. After this, we characterize graphs which are the exact-distance square of some bipartite graph and graphs which are the exact-distance square of some triangle-free graph.
	
	Our reduction is from the clique-edge-cover problem defined as follows. Given a graph $G$, a collection $C_{1},C_{2},\ldots,C_{k}$ of its clique subgraphs is said to form \emph{$k$-clique edge cover} if each edge of $G$ belongs to at least one $C_{i}$, $1\leq i \leq k$. The  \emph{clique-edge-cover problem} takes as input a pair $(G, k)$  of a graph $G$ and a positive integer $k$ and outputs YES if $G$ admits a $k$-clique edge cover, NO otherwise. 
	It is shown in \cite{Kou1978Covering} and \cite{Olrin1977} that the clique-edge-cover problem is an NP-complete problem. Here we show a polynomial time reduction from this problem to the $\mathcal{B}$-exact-distance-root problem, proving that the problem of deciding if a given graph admits a bipartite exact-distance square root is also NP-hard.
		
	\begin{proof}[Proof of Theorem~\ref{theorem:bipartite}]
		Note that the problem is in NP since computing exact-distance square of a given (bipartite) graph can be done in polynomial time. 
		To complete the proof, given a pair $(G, k)$ of a connected graph and a positive integer $k$, we build an auxiliary graph $G_k$ whose order is polynomial in the order of $G$ and $k$ and we show that $G$ admits a $k$-clique edge cover if and only if~$G_k$ admits bipartite exact-distance square root. We note that the assumption on connectivity does not change the nature of the $k$-clique edge cover problem, because if $G$ consists of two disconnected parts $G_1$ and $G_2$, we only need to consider the problems $(G_1, k_1)$ and $(G_2, k-k_1)$ for all choices of $0\leq k_1\leq k$.  
		
		The graph $G_k$ is built from $G$ by first adding a universal vertex to it, i.e, a new vertex $u$ which is joined to all vertices in $G$, and then adding a disjoint $k$-clique whose vertices we label $c_1, c_2, \ldots, c_k$. Clearly the order of $G_k$ is polynomial in terms of the input $(G, k)$. What remains to show is that the graph $G_k$ admits a bipartite exact-distance square	root $B$ if and only if $G$ admits a $k$-clique edge cover.
		
		To prove this claim we first assume that $G$ has a $k$-clique edge cover $C_{1},\ldots,C_{k}$. An example of a bipartite exact-distance square root $B$ of $G_k$ is as follows. Vertices of $B$ are $V(G)\cup \{u,c_{1},c_{2},\ldots,c_{k}\}$ with $V(G)\cup \{u\}$ forming one part and $\{c_{1},c_{2},\ldots,c_{k}\}$ forming the other part. Edges of $B$ consists of $E_{B}=\{uc_{i}:1\leq i\leq k\}\cup \{c_{i}v:v\in C_{i},1\leq i\leq k\}$. Let us verify that $G_k=\exact{B}{2}$. An adjacent pair $(x,y)$ of vertices of $G$ must be in clique $C_j$, thus they are both adjacent to the vertex $c_j$ of $B$ and hence adjacent in $\exact{B}{2}$. A nonadjacent pair $(x,y)$ of vertices of $G$ do not belong to any clique, and thus have no common neighbor in $B$, consequently they are not adjacent in $\exact{B}{2}$. As each vertex $c_i$ is adjacent to $u$ (in $B$), they form a clique in $\exact{B}{2}$. Finally, since $G$ is assumed to be connected, each vertex $x$ of $G$ is incident to some edge and thus belong to some $C_j$. That implies that $xc_ju$ is path of length 2 in $B$, concluding that $u$ is adjacent to all vertices $x$ of $G$ in $\exact{B}{2}$. In conclusion we have $G_k=\exact{B}{2}$.
		
		To complete the proof, we assume that there exists a bipartite graph $B$ satisfying  $G_k=\exact{B}{2}$. As $G$ is connected and, in $G_k$, the vertex $u$ is connected to vertices of $G$, vertices in $V(G)\cup \{u\}$ must all be in the same part of $B$. On the other hand, each pair $(c_i,c_j)$ of vertices in $\{c_1, c_2, \ldots, c_k\}$ must have a common neighbor which cannot be in $\{c_1, c_2, \ldots, c_k\}$ as it induces a clique in $\exact{B}{2}$. Thus each such pair has a common neighbor in part $V(G)\cup \{u\}$ and hence they form the other part of $B$. Let $C_i=N_{_B}(c_i)-u$. We claim that $\{C_1, C_2, \ldots, C_k\}$ forms a $k$-clique edge cover of $G$. That $C_i$ induces a clique in $G_k$ is consequence of the fact that its elements are vertices from same part of $B$, thus no two of them are adjacent, but as any two are adjacent to $c_i$, they are pairwise at distance 2 in $B$. That any edge $xy$ is in at least one $C_i$ is by the fact that they must be at distance 2 in $B$ and their common neighbor can only be a vertex in $\{c_1, c_2, \ldots, c_k\}$. 
	\end{proof}

	Though we have just shown that deciding if a graph has a bipartite exact-distance square root is an NP-complete problem, we  provide a characterization of these graphs as it may help in the development of new results or algorithms. This characterization is first guided by the observation that in the exact-distance square of a bipartite graph $B$ there is no connection between vertices in different parts, and that the neighborhood of each vertex in $B$ forms a clique in $\exact{B}{2}$. The exact statement of the characterization, Theorem~\ref{thm:Clique-dualCharacterization}, is based on the notion of clique-dual pairs defined as follows. Given a pair $F$ and $F'$ of graphs with no isolated vertices,  we say $F'$ is a \emph{clique-dual} of~$F$ if it admits a clique edge cover $Cl_{1},\ldots,Cl_{n}$ labeled by the vertices of $F$ such that each edge of~$F'$ belongs to at least one clique and $v_i\sim v_j$ in $F$ if and only if $Cl_{i}\cap Cl_{j}\neq \emptyset$.
	It is not hard to build examples of nonisomorphic clique-duals or graphs which admit no clique-dual. 
	However, the following holds.
	
	\begin{proposition}
		Given graphs $F$ and $F'$ if $F'$ is a clique-dual of $F$, then $F$ is also a clique-dual of~$F'$.
	\end{proposition}	 
	
	\begin{proof}
		Let $v_1, v_2, \ldots, v_n$ be the vertices of $F$ and $Cl_{1},\ldots,Cl_{n}$ a set of cliques in $F'$ presenting it as a clique-dual of $F$. Let $u_1, u_2, \ldots, u_k$ be the vertices of $F'$. Given a vertex $u_i$ of $F$, let $Cl_{i_1},\ldots,Cl_{i_l}$ be the set of cliques each of which contains $u_i$. We claim that the set $\{v_{i_1}, v_{i_2}, \ldots, v_{i_l}\}$ of the vertices in $F$ form a clique $Cl'_{i}$. That is because the cliques corresponding to any two of these vertices have at least $u_i$ in common. Let $v_av_b$ be an edge of $F$. Then there exist a $u_r\in Cl_a \cap Cl_b$. Hence $Cl'_r$ contains both $v_a$ and $v_b$. Thus the collection $Cl'_{i}$, $i=1, 2,\ldots, k$, is a clique edge cover of $F$. It remains to show that vertices $u_i$ and $u_j$ of $F'$ are adjacent (in $F'$) if and only if $Cl'_i\cap Cl'_j\neq \emptyset$. If $u_i\sim u_j$, then, since $Cl_{1},\ldots, Cl_{n}$ is an edge cover of $F'$, there is a clique, say $Cl_1$, that contains both $u_i$ and $u_j$. Then, by our choice, $Cl'_{i}$ and $Cl'_{j}$ both contain the vertex $v_1$. For the inverse, assume $v_1$ is a common element of $Cl'_{i}$ and $Cl'_{j}$. That means $u_i$ and $u_j$ are both in $Cl_1$, but then they must be adjacent. 			
	\end{proof}

	\begin{theorem}\label{thm:Clique-dualCharacterization}
		A graph $G$ with no isolated vertex	has a bipartite exact-distance square root if and only if it consists of two disjoint parts $F$ and $F'$ such that there is no connection between the two and that $F'$ is a clique-dual of $F$.
	\end{theorem}
	
	\begin{proof}
		Given a connected bipartite graph $H$ with $X$ and $Y$ being the vertex sets of the two parts, $\exact{H}{2}$ consists of two connected parts one induced by $X$ and another induced by $Y$. These two are easily observed to form a clique-dual pair. Thus the only if part of the theorem follows easily. 
		
		Now let $F$ and $F'$ be the pair of clique-dual graphs forming $G$. Assume that $v_1, v_2, \ldots, v_n$ are the vertices of $F$, and let $u_1, u_2, \ldots, u_k$ be the vertices of $F'$. Furthermore, let $Cl_{1},\ldots, Cl_{n}$ be the cliques of $F'$ providing the clique-dual relation. Let $H$ be the bipartite graph with $v_1, v_2, \ldots, v_n$ forming vertices of one part and $u_1, u_2, \ldots, u_k$ forming the vertices of the other part, with $v_i$ being adjacent to all $u_j$ in $Cl_{i}$. We claim that $G=\exact{H}{2}$. We first consider the subgraph of $G$ induced by vertices $v_1, v_2, \ldots, v_n$. If vertices $v_i$ and $v_j$ are adjacent in $G$, then the clique $Cl_{i}$ and $Cl_{j}$ must have a common vertex, say $u_{ij}$. Thus in $H$ both $v_i$ and $v_j$ are connected to $u_{ij}$ and, hence, noting that $H$ is bipartite, they are at distance 2 in $H$. Conversely, if $v_i$ and $v_j$ are not adjacent in $G$, then they cannot have a common neighbor in $H$ as otherwise they must belong to the corresponding clique of the common vertex. Thus $v_1, v_2, \ldots, v_n$ induce the component $F$ of $G$. Next we consider vertices $u_1, u_2, \ldots, u_k$. If $u_i$ is adjacent to $u_j$ in $G$, then they are in a clique $Cl_l$ of $F'$. Then they have $v_l$ as a common neighbor in $H$, and, as $H$ is bipartite, they are at distance 2 in $H$. If $u_i$ and $u_j$ are not adjacent in $G$, then they cannot have a common neighbor in $H$, as if a vertex, say $v_l$, was a common neighbor, then they both must have been in $Cl_l$, which means they must have been adjacent in $G$.	
	\end{proof}

	We end this section by providing a characterization of exact-distance squares of triangle-free graphs, which has a very similar taste to the one of the previous theorem.

	\begin{theorem}\label{thm:Triangle-free-Root}
		A graph $G$ on vertices $v_1,v_2, \ldots, v_n$ admits a triangle-free exact-distance square root if and only if it has a collection of cliques $Cl_{1}, Cl_{2} \ldots,Cl_{n}$ satisfying the following properties: 
		\begin{itemize}
			\item $v_i\not\in Cl_{i}$, for each $i$
			\item if  $v_i\in Cl_{j}$, then  $v_j\in Cl_{i}$ for each pair $i,j$,
			\item $v_i\sim v_j$ if and only if $Cl_{i}\cap Cl_{j}\neq \emptyset$.
			\item For each $i\in [n]$, and for each pair $v_j, v_k$ of vertices of $Cl_{i}$ we have $v_j\not \in Cl_{k}$.
		\end{itemize}
		
	\end{theorem}

	\begin{proof}
		For the easier direction, assume $H$ is a triangle-free exact-distance square root of $G$, thus, by the definition, $V(H)=V(G)$. It is then enough to take $Cl_{i}=N_{H}(v_i)$. Observe that if $H$ is triangle-free, then  $N_{H}(v_i)$ induces a clique in $G$. If $v_i\sim v_j$ (in $G$), then in $H$ they are at distance 2. Let $x$ be a common neighbor, then $x\in Cl_{i}\cap Cl_{j}$.  Conversely, if there is a vertex $x$ in $Cl_{i}\cap Cl_{j}$, then $v_i$ and $v_j$ are at distance at most 2 in $H$, but since $H$ is triangle-free, they are indeed at distance exactly 2 and thus adjacent in $G$. To check that the last condition holds, note that if $v_j$ and $v_k$ are both neighbors of $v_i$, since $H$ is triangle-free, $v_j$ is not adjacent to $v_k$ in $H$ and, therefore, $v_j\not \in Cl_{k}$.
		
		For the other direction, let $Cl_{1}, Cl_{2} \ldots,Cl_{n}$ be a collection of the cliques satisfying the properties. We construct a graph $H$ as follows and show that it is an exact-distance square root of $G$ and that it has no triangle. For the vertices we take $V(H)=V(G)$. For the edges, each vertex $v_i$ is adjacent to all vertices in $Cl_{i}$. The first condition implies that $H$ has no loop. The second condition implies that $Cl_{i}=N_{H}(v_i)$. The third condition implies that if $v_i$ and $v_j$ are at distance 2 in $H$, then they are adjacent in $G$ because their common neighbor(s) will be in both $Cl_i$ and $Cl_j$, thus indeed $\exact{H}{2}=G$. Finally, using the last condition we show that $H$ must be triangle-free. That is because, if $v_i$ is adjacent to vertices $v_j$ and $v_k$ in $H$ (i.e., $v_j, v_k \in Cl_i$), then by the last condition $v_j\not \in Cl_k$ and, therefore, by the construction of $H$, $v_j$ is not adjacent to $v_k$ in $H$.  
	\end{proof}

	\section{Characterizing exact-distance squares of trees}\label{sec:charactrees}
	
	In this section we prove Corollary~\ref{corollary:trees}, a characterization of exact-distance squares of trees. We will use this characterization, in the next section, to provide  a polynomial-time algorithm to recognize exact-distance squares of trees.

	\begin{theorem}
		\label{twoCliqueTrees}
		Let $T$ be a non-trivial tree. Then $\exact{T}{2}$ consists of two connected components each of which is a clique-tree.
	\end{theorem}
	
	\begin{proof}
		Since $T$ is connected and bipartite, $\exact{T}{2}$ consists of two connected components, say $T_1$ and $T_2$, one induced by each part of $T$. Applying the proof of Theorem~\ref{thm:Clique-dualCharacterization} $T_1$ is covered by a set of cliques that are labeled by the vertices of $T_2$, where the clique corresponding to the vertex $v$ is induced by the neighbors of $v$ in $T$. 
		
		To see that $T_1$ is a clique-tree, assume (toward a contradiction) that there is cycle in $T_1$ which does not induce a complete graph and let $C$ be a shortest of all such cycles. Then, on the one hand, $C$ uses at most two vertices from each clique of $T_1$ as otherwise one can find a shorter cycle, noting moreover that any such a pair of vertices must be consecutive vertices of $C$. On the other hand, each edge $uv$ of $C$ must be in one of the selected cliques because it represents a pair of common nieghbors of a vertex, say $x$ in $T_2$. Hence, $x$ is not adjacent to any other vertex of $C$. Then we consider the cycle $C'$ obtained from $C$ where each edge $uv$ is replaced with the 2-path $u-x-v$. This cycle $C'$ however must be a cycle of $T$, contradicting the fact that $T$ is a tree.     
	\end{proof}
	
	So among all graphs which are disjoint union of two clique-trees, we aim at distinguishing which are the exact-distance squares of trees. cut-vertices will naturally play a key role, and in particular we will need the following.

	\begin{lemma}
		\label{cutvertex}
		A vertex $v$ of tree $T$ is a cut-vertex of $\exact{T}{2}$ if and only if it has at least two non-leaf neighbors in $T$.   
		
	\end{lemma}
	\begin{proof}
		As $\exact{T}{2}$ is a union of two clique-trees, if $v$ is a cut-vertex of $\exact{T}{2}$, then there are at least two blocks $B_1$ and $B_2$ such that $v\in B_1\cap B_2$. Each of these two blocks then corresponds to the set of neighbors of vertices say $w_1$ and $w_2$. These are the two neighbors of $v$ that are not leaves.
		Conversely, if $w_1$ and $w_2$ are two non-leaf neighbors of $v$, then the cliques $B_1$ and $B_2$ corresponding to $w_1$ and $w_2$ will become disconnected after removing $v$.
	\end{proof}

	\begin{corollary}
		\label{degreeVSblocks}
		If $v\in\exact{T}{2}$ is in $k$ blocks of $\exact{T}{2}$, then $v$ has at least $k$ neighbours in $T$.
	\end{corollary}

	Let $W$ be a clique-tree. We define the \emph{canonical tree} of $W$, denoted, $T_W$ as follows. For each block $B$ in $W$ we create a new vertex $v_B$ and set $V(T_W)=V(W)\cup\{v_B|B\text{ is a block of }W\}$, and $E(T_W)=\{\{u,v_B\}|B\text{ is a block of }W\text{ and }u\in B\}$. If $v_B$ is the vertex in $T_W$ corresponding to the block $B$ of $W$, then we say that $v_B$ \emph{arises from $B$ in $T_W$}. We first show that these graphs are indeed trees.
	
	\begin{lemma}
		Let $W$ be a clique-tree. Then $T_W$ is a tree.
	\end{lemma}
	
	\begin{proof}
		Let $\{u_1,\ldots,u_n\}$ be the vertices of $W$ and $\{B_1,\ldots,B_m\}$ its blocks. Since $W$ is connected it is not hard to see that $T_W$ is connected. Now, towards a contradiction, we assume that $T_W$ has a cycle. Let $C$ be a shortest cycle in $T_{W}$. The cyclic ordering of the vertices of $C$ alternates in the following form: $u_{i_1}v_{B_{j_1}}u_{i_2}v_{B_{j_2}}\ldots u_{i_k}v_{B_{j_k}}u_{i_1}$. It follows that $u_{i_1}u_{i_2}\ldots u_{i_k}u_{i_1}$ is cycle of $W$, but then they are in the same block, contradicting that $C$ is a cycle of $T_W$. 
	\end{proof}
	
	We pause briefly to reap the following characterization of clique-trees.
	
	\begin{theorem}
		A graph is a clique-tree if and only if it is isomorphic to a component of the exact-distance square of some tree.
	\end{theorem}
	\begin{proof}
		It is not hard to see that every clique-tree $W$ is a component of $\exact{T_W}{2}$ (and, by the previous lemma, $T_W$ is a tree). The other direction is given by Theorem~\ref{twoCliqueTrees}.
	\end{proof}
	
\begin{figure}[ht]
		\centering
		\includegraphics[height=3.2 in]{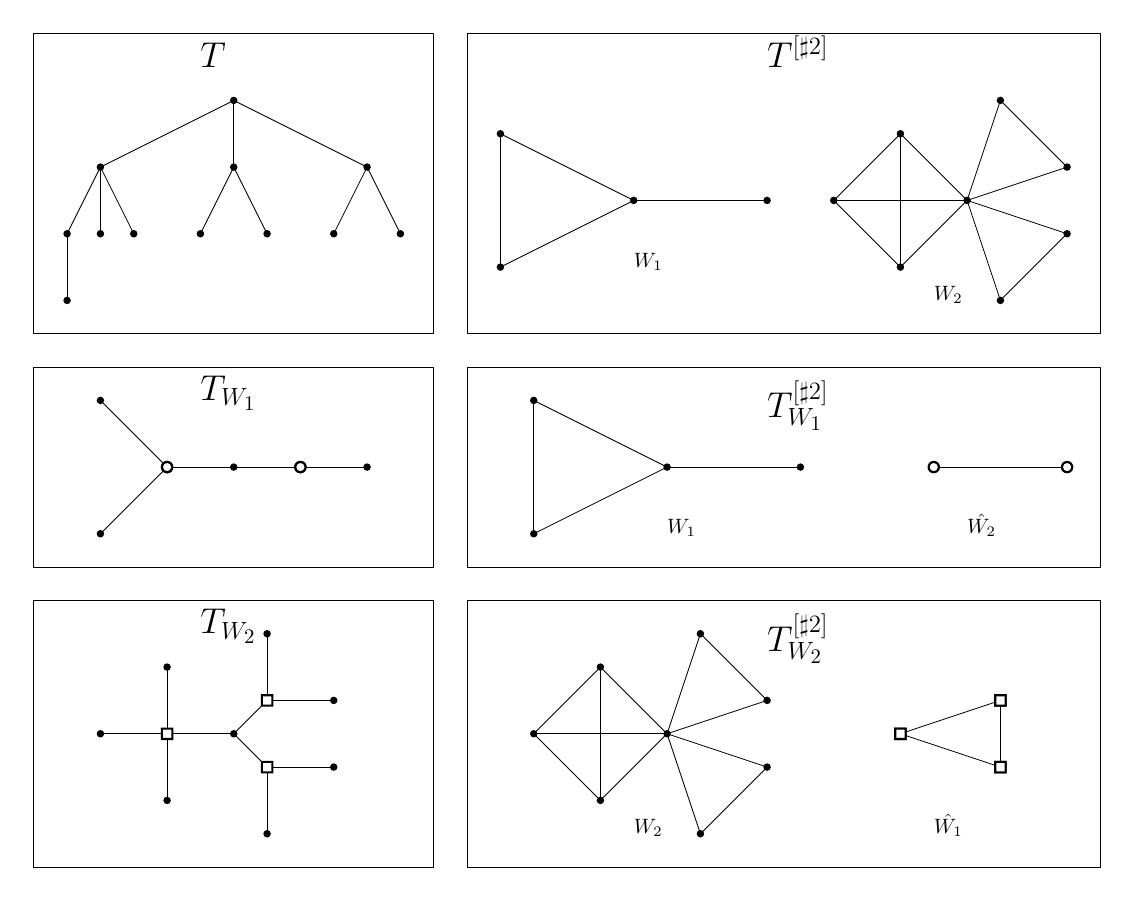}
		\caption{A tree $T$ and its exact-distance square with components $W_1$ and $W_2$. Here $W_1$ is not isomorphic to any component of $\exact{T_{W_{2}}}{2}$, and $W_2$ is not isomorphic to any component of $\exact{T_{W_{1}}}{2}$.}
		\label{I1}
	\end{figure}

	We now want to tell whether the disjoint union of two given clique-trees $W_1$ and $W_2$ is the exact-distance square of some tree. We would like to assign to each clique of $W_1$ a vertex $v$ in $W_2$ such that the clique would represent the neighborhood of $v$ in a tree. This, of course, is reminiscent of the way we construct the canonical tree of a clique-tree. Indeed, if we consider the exact-distance square of $T_{W_1}$, one of its components will be $W_1$, and if the other one turns out to be isomorphic to $W_2$, then we can immediately tell that $W_1$ and $W_2$ form the exact-distance square of some tree. But it might happen that $W_1$ and $W_2$ form the exact-distance square of some tree but $W_1$ is not isomorphic to any component of $\exact{T_{W_{2}}}{2}$, and $W_2$ is not isomorphic to any component of $\exact{T_{W_{1}}}{2}$ (see Figure~\ref{I1}, for an example). In this case, however, the component of $\exact{T_{W_{1}}}{2}$ which is not $W_1$ will always be a subgraph of $W_2$. In such a case if we find an adequate isomorphism from this component to some subgraph of $W_2$, then we will be able to ``complete" (through Procedure 1) the canonical tree of $W_1$ into a tree exact-distance square root of the union of $W_1$ and $W_2$. What we mean by an adequate isomorphism is made explicit in the theorem below, which is the last step towards our characterization of exact-distance squares of trees. Here, given a subset $S$ of the vertices of graph~$G$, the set of cut-vertices of $G$ in $S$ is denoted by $\CV (S)$.

	\begin{theorem}\label{nearcharac}
		Let $G$ be a disjoint union of two clique-trees $W_1$ and $W_2$. Let $\hat{W}_2$ be the component of $\exact{T_{W_1}}{2}$ which is not $W_1$. Then there exists a tree $T$ such that $\exact{T}{2}=G$ if and only if there is an isomorphism $\varphi$ of $\hat{W}_2$ to a subgraph of $W_2$ satisfying the following two conditions:

		\begin{itemize}
			\item[(1)] For every vertex $v\in\hat{W_2}$ we have that $d_{T_{W_1}}(v)$ is at least the number of blocks to which $\varphi(v)$ belongs in $W_2$
			\item[(2)] If we have $x\in \CV(W_2)$, then we have $x\in\varphi(\hat{W_2})$
		\end{itemize}
	\end{theorem}  
	
	A key element of our proof of this theorem is Procedure~1, which given an isomorphism satisfying conditions (1) and (2) ``completes" $T_{W_1}$ into a tree $T$ having $G$ as its exact-distance square. Before going into the proof let us exemplify how this algorithm works. In Figure~\ref{I2} we have a graph $G$ which is the disjoint union of two clique-trees $W_1$ and $W_2$. Let us assume that this graph is given as an input for Procedure 1, together with $T_{W_1}$ and an isomorphism $\varphi$ satisfying the conditions of the theorem. The component $W_2$ has exactly one cutvertex called $x$. The preimage of $x$ through $\varphi$ is $v$. (Note that $w$ could not be the preimage of $x$, because otherwise Condition (1) of the theorem would not be satisfied). In the proof, for a block $B$ we let $A_B$ be the set of the vertices in $B$ which have no preimage by $\varphi$.  In $W_2$ the only block that is considered in the second {\bf for} loop is the block that contains $y$ (for $y\in B\setminus(A_B\cup\{x\})\neq\varnothing$ if $B$ is that block). This loop creates a new tree $T$ from $T_{W_1}$  by adding each vertex of the block of $y$ that does not have a preimage, (in this case just $a_1$) and making it adjacent to $z$. The other two blocks of $W_2$ are considered in two subsequent iterations of the third {\bf for} loop, and the image shows how the tree $T$ is updated in this iterations. Note that the final tree obtained satisfies $\exact{T}{2}=G$.

	\begin{proof}
		We first assume that there exists $T$ such that $\exact{T}{2}=G$. Consider an arbitrary vertex $v\in V(\hat{W_2})$. This vertex can be seen as arising in $T_{W_1}$ from a block $B_v$ of $W_1$. 
		
		Since $G$ is exact distance square of a bipartite graph ($T$), and following Theorem~\ref{thm:Clique-dualCharacterization}, $W_1$ and $W_2$ are clique-dual. Thus there exists a vertex $x$ in $W_2$ such that $B_v=W_1[N_{T}(x)]$, noting that this set is uniquely determined because $T$ is a tree. Let $$S =\{x\in V(W_2) \mid \exists v\in\hat{W_2}\text{ which arises from }W_1[N_{T}(x)]\text{ in }T_{W_1}\}.$$ We define $\varphi:V(\hat{W_2})\rightarrow S$ such that $$\varphi(v)=x \mbox{ if and only if } v  \mbox{ arises from } W_1[N_{T}(x)]\mbox{ in }T_{W_1}.$$ Since a block is induced by the neighborhood of a unique vertex, $\varphi$ is well defined and is a bijection.
		
		We first show that $\varphi$ is an isomorphism. Having $uv\in E(\hat{W_2})$ is equivalent to having two blocks $B_u$ and $B_v$ in $W_1$, such that $u$ and $v$ arise in $T_{W_1}$ from $B_u$ and $B_v$, respectively, and satisfy $B_u\cap B_v\neq\varnothing$. But $B_u$ and $B_v$ are induced in $W_1$ by the neighborhood in $T$ of $\varphi(u)$ and $\varphi(v)$, respectively, and if these two neighborhoods intersect, then we have $\varphi(u)\varphi(v)\in E(W_2)$. So $uv\in E(\hat{W_2})$ is equivalent to $\varphi(u)\varphi(v)\in E(W_2)$.
		
		To show we have $(1)$, let $v\in V(\hat{W_2})$. Notice that $d_{T_{W_1}}(v)=|V(B_v)|$, where $B_v$ is the block in $W_1$ such that $v$ arises from $B_v$ in $T_{W_1}$. By definition of $\varphi$, $B_v=W_1[N_T(\varphi(v))]$. This easily gives $d_T(\varphi(v))=|V(\hat{B_v})|$. Since we have $d_{T_{W_1}}(v)=d_T(\varphi(v))$, by Corollary \ref{degreeVSblocks} we obtain $(1).$ To prove (2) we take a vertex $x\in \CV (W_2)$. Notice that the neighborhood of $x$ in $T$ induces a block $B$ in $W_1$, and then there exists a vertex $v_{B}\in \hat{W_2}$ which arises from $B$ in $T_{W_1}$. Then, by definition of $\varphi$, we have that $\varphi(v_B)=x$.
		
		\begin{figure}[ht]
			\centering
			\includegraphics[height=12cm]{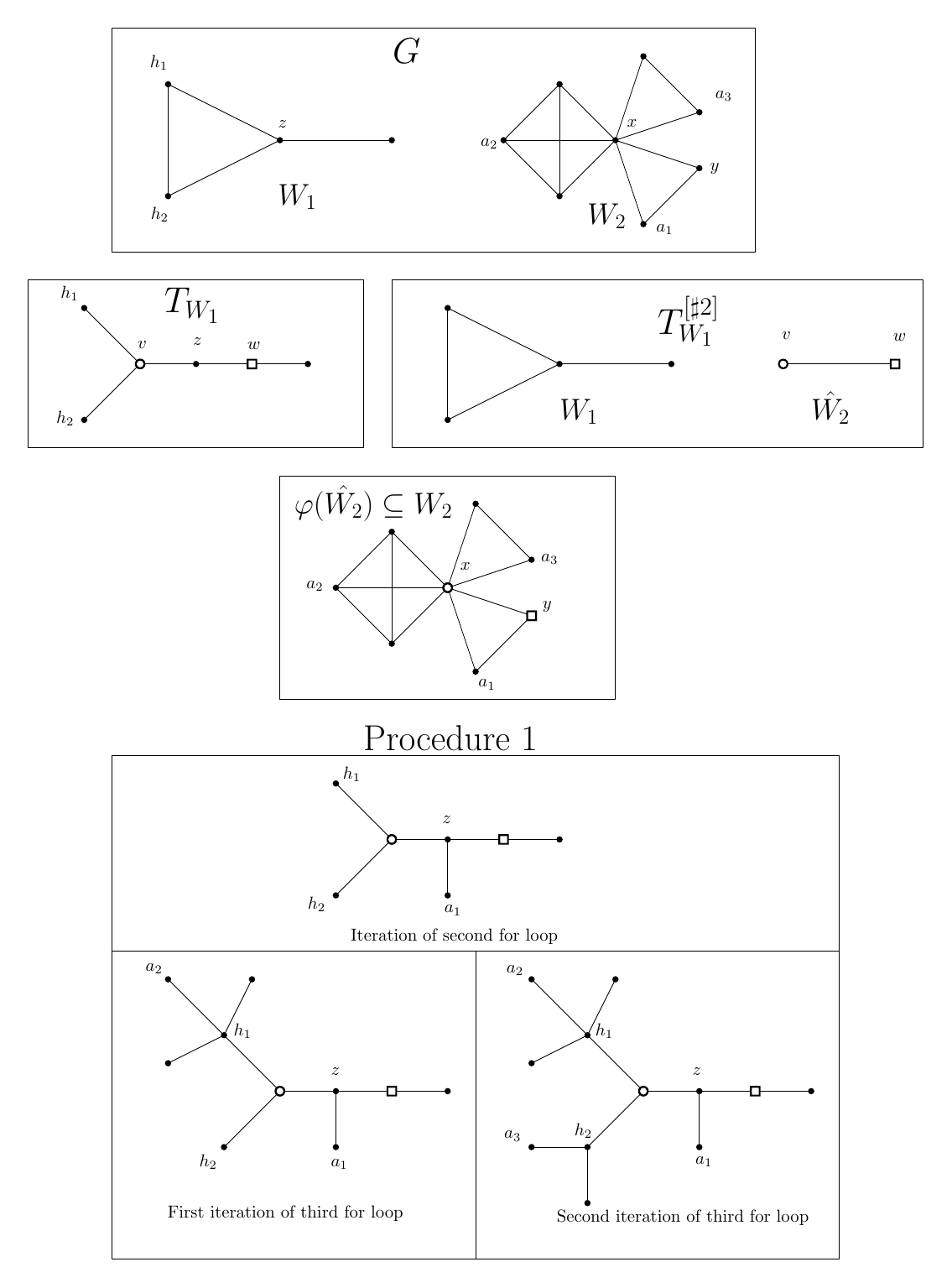}
			\caption{Constructing, through Procedure 1, a tree which is the  exact-distance square root of a given graph $G$}
			\label{I2}
		\end{figure}
		
		For the other direction, we assume that there exists an isomorphism $\varphi$ of $\hat{W_2}$ to a subgraph of $V(W_2)$ which satisfies $(1)$ and $(2)$. Then starting from $T_{W_1}$, we construct a tree $T$ such that $\exact{T}{2}=G$.
		
		Notice that since $\varphi$ is an isomorphism, and $W_2$ and $\hat W_2$ are clique-trees, we have $\varphi(\CV (\hat{W_2}))\subseteq \CV (W_2)$. First we consider the case where $W_2$ has no cut-vertex. In this case $\hat{W_2}$ has no cut-vertex either. Then $W_2$ and $\hat{W_2}$ have only one block and $|V(\hat{W_2})|\leq|V(W_2)|$. Let $w\in W_1$ be the unique vertex whose neighborhood in $T_{W_1}$ induces $\hat{W_2}$. We add to $T_{W_1}$ all the vertices of $W_2\setminus\varphi(\hat{W_2})$ as leaf-neighbours of $w$. The new tree satisfies that its exact-distance square is isomorphic to~$G$.
		
		Suppose then that $W_2$ has at least one cut-vertex. By (2) we have that for any $x\in \CV (W_2)$ there exists $v\in\hat{W_2}$ such that $\varphi(v)=x$. For every block $B$ of $W_2$ we let $A_B=\{w\in B|w\in W_2\setminus\varphi(\hat{W_2})\}$, that is, $A_B$ is the set of the vertices in $B$ each of which does not have a preimage by $\varphi$.
		
		Based on the notation, we now introduce Procedure 1 using which we can build the tree $T$ whose exact distance square is $G$. The proof of the correctness of the procedure is followed.
		
		\begin{algorithm}[ht]
			\begin{algorithmic}[1]
				
				\Require $T_{W_1}, G, \varphi$ as in Theorem~\ref{nearcharac}
				\Ensure $T$ tree such that $\exact{T}{2}\cong G$
				\State $F \gets \varnothing$
				\State $T \gets T_{W_1}$
				\For{$x\in \CV(W_2)$}
				\State $v \gets \varphi^{-1}(x)$
				\For{$B$ block in $W_2$ such that $x\in B$ and $B\setminus(A_B\cup\{x\})\neq\varnothing$}
				\State Let $y\in B\setminus(A_B\cup\{x\})$ 
				\State $w\gets \varphi^{-1}(y)$
				\State Let $z\in T$ be such that $z\in N_{T_{W_1}}(v)\cap N_{T_{W_1}}(w)$
				\If{$z\notin F$}
				\State $V(T)\gets V(T)\cup A_B$
				\State $E(T)\gets E(T)\cup \{az|a\in A_B\}$
				\State $F \gets F\cup\{z\}$
				\EndIf
				\EndFor
				\For{$B$ block in $W_2$ such that $x\in B$ and $B\setminus(A_B\cup\{x\})=\varnothing$}
				\State Let $h\in N_{T_{W_1}}(v)\setminus F$ 
				\State $V(T)\gets V(T)\cup A_B$
				\State $E(T)\gets E(T)\cup \{ah|a\in A_B\}$
				\State $F \gets F\cup\{h\}$
				\EndFor
				\EndFor
				\State \Return $T$
			\end{algorithmic}
			\caption{EXACT-DISTANCE SQUARE ROOT TREE COMPLETION}
		\end{algorithm} 
		
		Consider $x\in \CV (W_2)$ and let $v$ be its preimage by $\varphi$. In the case where $B$ satisfies $B\setminus(A_B\cup\{x\})\neq\varnothing$ we take $y\in B\setminus(A_B\cup\{x\})$ with its preimage~$w$. Since $\varphi$ is an isomorphism we have that $vw\in E(\hat{W_2})$ and then we can take the unique vertex $z\in N_{T_{W_1}}(v)\cap N_{T_{W_1}}(w)$. If $z\notin F$, then we add the vertices of $A_B$ to $T$, each as a leaf-neighbor of $z$. If $T$ was already a tree, then this operation clearly keeps it so: indeed, since we must have $z\in V(W_1)$, in this case we only add leaves to a vertex of $T_{W_1}$.

		Observe that for the chosen vertex $x$, and for each block that contains $x$, the procedure will add at most one vertex to $F$. Then, by condition (1) of the theorem, there is always a vertex $h\in N_{T_{W_1}}(v)\setminus F$ in the case when $B\setminus(A_B\cup\{x\})=\varnothing$. To this vertex $h$ we join the vertices of $A_B$ as leaves. Since we have $x\in \CV(W_2)\subseteq V(W_2)$ and $h\in N_{T_{W_1}}(v)$, we must have $h\in V(W_1)$. Thus, again, all we do in the third {\bf for} loop is add leaves to a vertex of $T_{W_1}$. Altogether, the output $T$ is indeed a tree. Moreover, since we assume $\CV (W_2)\neq\varnothing$, we have that every block in $W_2$ contains at least one cut-vertex, so in this process we add to $T_{W_1}$ all the vertices in $W_2\setminus\varphi(\hat{W_2})$ as leaves.

		To finish the proof, we need to show that $\exact{T}{2}$ is isomorphic to $G$. We define $\psi:V(\exact{T}{2})\rightarrow V(G)$ such that $\psi(v)=v$ if $v\in V(W_1)\cup V(W_2)\setminus(V(\varphi(\hat{W_2})))$, or $\psi(v)=\varphi(v)$ if $v\in V(\hat{W_2})$. Let us see that $\psi$ is an isomorphism. It is clear that $\psi$ is a bijection. Now let $\psi(u)\psi(v)\in E(G)$. There are only four possible cases:
		
		Case 1: $\psi(u),\psi(v)\in V(W_1)$. Since $T$ was obtained from $T_{W_1}$ by merely adding leaves to vertices in $W_1$, having $\varphi(u)\varphi(v)=uv\in E(W_1)$ is equivalent to $uv\in E(\exact{T}{2})$.
		
		Case 2: $\psi(u),\psi(v)\in V(W_2)\setminus V(\varphi(\hat{W_2}))$. In this case $\varphi(u)\varphi(v)=uv\in E(W_2\setminus\varphi(\hat{W_2}))$ which is equivalent to that there exists a block $B$ in $W_2$ such that $u,v\in A_B$. By Procedure 1, this is equivalent to $uv\in E(\exact{T}{2})$.
		
		Case 3: $\psi(u),\psi(v)\in V(\varphi(\hat{W_2}))$. In this case we have $\psi(u)\psi(v)=\varphi(u)\varphi(v)\in E(\varphi(\hat{W_2}))$ and, since $\varphi$ is an isomorphism, this is equivalent to having $uv$ as an edge in $\hat{W_2}$. In turn, this is equivalent to having $uv\in E(\exact{T}{2})$.
		
		Case 4: $\psi(u)\in V(W_2)\setminus V(\varphi(\hat{W_2}))$ and $\psi(v)\in V(\varphi(\hat{W_2}))$. In this case $\psi(u)\psi(v)=u\varphi(v)\in E(W_2)$. This is equivalent to there being a block $B$ in $W_2$ containing $u\varphi(v)$ and such that $u\in A_B$. But in the algorithm $u$ is added to $T$ as a neighbour of a vertex $z$ which is neighbour of $v$. In other words, $u\varphi(v)\in E(W_2)$ is equivalent to $uv\in E(\exact{T}{2})$.
		
		The result follows.
	\end{proof}

	\begin{corollary}
		\label{fullcharac}\label{corollary:trees}
		A graph is an exact-distance square of some non-trivial tree if and only if it has two components $W_1$ and $W_2$, each of them a clique-tree, such that there exists an isomorphism $\varphi$ from $\exact{T_{W_1}}{2}\setminus W_1$ to a subgraph of $V(W_2)$ such that 
		\begin{itemize}
			\item[(1)] For every vertex $v\in V(T_{W_1}\setminus W_1)$ we have that $d_{T_{W_1}}(v)$ is at least the number of blocks to which $\varphi(v)$ belongs in $W_2$.
			\item[(2)] If we have $x\in \CV(W_2)$, then  $x$ belongs to the image of $\varphi$.
		\end{itemize}
	\end{corollary}

	\section{Recognizing exact-distance squares of trees}\label{sec:recognition}

	In this section we use the last stated result to prove Theorem~\ref{theorem:trees}, that is, to give a polynomial-time algorithm to decide, for a given a graph $G$, if there exists a tree $T$ such that $\exact{T}{2}=G$.

	In \cite{matula1978subtree} Matula gave a polynomial-time algorithm to decide if a tree $S$ on $n_S$ vertices is isomorphic to any subtree of a tree $T$ on $n_T$ vertices. We modify  ``Algorithm A'' of his paper to obtain an algorithm which, given a graph $G$, decides if $G$ has an exact-distance square root which is a tree. We note at the end how this algorithm may also produce one such a tree.
	
	We will first sketch Matula's algorithm and for this we need some definitions and a key result from \cite{matula1978subtree}. For a tree $T$ and $r\in V(T)$, we denote by $T[r]$ the tree $T$ rooted at $r$. For $u,v$ such that $uv\in E(T)$, we define the \emph{limb} $T[u, v]$ as the maximal subtree of $T$ containing the edge $uv$ in such a way that $u$ is a leaf of the subtree, and where $u$ is assigned as the root. The \emph{height} of a limb $T[u,v]$ is the largest distance of any vertex of $T[u,v]$ to its root $u$. For a rooted tree $T[u]$, the \emph{limbs of $T[u]$} are the limbs $T[v, w]$ where $v$ is either $u$ or in between $u$ and $w$.
	
	Let $S$ and $T$ be trees with $x,y\in V(S)$ and $u,v\in V(T)$. The rooted tree $S[x]$ is isomorphic to the rooted tree $T[u]$ if there is an isomorphism of $S$ to $T$ which takes $x$ into $u$. An isomorphism of the limb $S[x, y]$ to the limb $T[u, v]$ is called a \emph{limb embedding} of $S[x, y]$ in $T[u, v]$.  If there is such a limb embedding, then we say that $S[x, y]$ \emph{can be embedded} in $T[u, v]$. Let $a_1, a_2,\ldots,a_p$ be the neighbors of $y$ in $S[x, y]-x$, and $b_1,b_2,\ldots,b_q$, the neighbors of $v$ in $T[u, v]-u$. The \emph{highest-limbs matrix} associated with $S[x, y]$ and $T[u, v]$ has rows corresponding to the limbs $S[y, a_i]$, $1\leq i\leq p$, and columns corresponding to the limbs $T[v, b_j]$, $1\leq j\leq q$,	and its  position $(S[y, a_i], T[v, b_j])$ has value 1 if $S[y, a_i]$ can be embedded in $T[v, b_j]$ and 0 otherwise.
	
	For a  $0,1$-matrix on $p\times q$ entries, $p\leq q$, a \emph{bipartite matching} is a set of unit entries no two of which are from the same column or from the same row. A maximum bipartite matching is \emph{complete for the rows} if there are $p$ unit entry positions in the bipartite matching, one for each row.
	
	The following theorem is central to the subtree-isomorphism algorithm in \cite{matula1978subtree}.
	
	\begin{theorem}[Matula \cite{matula1978subtree}]
		\label{MatulaTheorem}
		The limb $S[x,y]$ is isomorphic to a rooted subtree of the limb $T[u,v]$ (with $u$ as the root) if and only if the associated highest-limbs matrix has a maximum bipartite matching which is complete for the rows.
	\end{theorem}
	
	To complete the ingredients of Matula's algorithm we need the following definition. The \emph{limb embedding matrix} $M(S[x,y], T)$ has a row for each of the $n_{S}-1$ limbs of the rooted tree $S[x, y]$ and a columns for each of the $2(n_{T}-1)$ limbs of the unrooted tree $T$. More specifically,  we order the limbs of $S[x, y]$ by non-decreasing height and order the rows accordingly. For all limb pairs $S[a, b]$ of $S[x, y]$ and $T[u, v]$ of $T$, the entry for the $(S[a, b], T[u, v])$ position is 1 if $S[a, b]$ can be embedded in $T[u, v]$ and 0 otherwise. 
	
	The subtree-isomorphism algorithm from \cite{matula1978subtree} works as follows. Note that deciding if there is a subtree of $T$ isomorphic to  $S$ is equivalent to deciding whether the limb $S[l,y]$ rooted at a leaf $l$ of $S$ can be embedded in some limb of $T$. It is not hard to see that the limb $S[l,y]$ can be embedded in some limb of $T$ if and only if there is an entry of 1 in the last row of the corresponding limb embedding matrix. So we need an efficient way of computing this matrix. For all rows of this matrix corresponding to height-1 limbs of $S[l,y]$, set all entries to 1. Now assuming that the matrix has been computed correctly up until some row, we build the next row, corresponding to a limb, say $S'$ of $S[l,y]$, as follows. For each entry we build the highest-limbs matrix of $S'$ and $T'$, where $T'$ is the limb of $T$ corresponding to the column of that entry. This matrix only needs information of limbs of smaller height, and thus we can obtain it from the previous rows of the limb embedding matrix. Deciding the value of the entry is then reduced, by Theorem~\ref{MatulaTheorem}, to a bipartite matching problem (recall that this can be solved in polynomial-time \cite{HopcroftTarjan}). Since we only need to solve a quadratic amount of these problems (in terms of max$\{|V(S)|, |V(T)|\}$), the algorithm runs in polynomial time. In fact, an isomorphism, if it exists, can be recovered in polynomial time, as we shall later mention.

	We now state our algorithm, which can be seen as having two stages. Its input is a graph $G$, and it decides if there exists some tree $T$ such that $\exact{T}{2}=G$.
	
	\vskip .1 in

	\noindent\textbf{Stage I. Initialization}
	
	In this stage we compute some relevant graphs, as well as important information we will use in the second stage.
	
	\begin{itemize}
		\item[I.1] Decide if $G$ has exactly two components $C_1$, $C_2$, each a clique-tree. If it does not, then return that there is no tree having $G$ as its	exact distance square.
		\item[I.2] For every vertex $v$ in $C_2$ we store the number $b_{C_2}(v)$ of blocks to which it belongs in $C_2$.
		\item[I.3] From $C_1$ and $C_2$ we construct $T_{C_1}$ and $T_{C_2}$.
		\item[I.4] From $T_{C_1}$ we construct $\exact{T_{C_1}}{2}$. Notice that $\exact{T_{C_1}}{2}$ consists of $C_1$ and another clique-tree $\hat{C_2}$. For every vertex $v$ in $\hat{C_2}$ we store $d_{T_{C_1}}(v)$.
		\item[I.5] From $\hat{C_2}$ we construct $T_{\hat{C_2}}$.
	\end{itemize}

	\noindent\textbf{Stage II. Obtaining the isomorphism}
	
	In this stage we determine if there is an isomorphism from $T_{\hat{C_2}}$ to a subtree of $T_{C_2}$ from which we can construct an isomorphism between $\hat{C_2}$ and a subgraph of $C_2$ which satisfies the conditions of Theorem \ref{nearcharac}.
	
	\begin{itemize}
		
		\item[II.1] Consider $T_{\hat{C_2}}$ as rooted at a leaf $l$ and let $y$ be its (only) neighbor. Order the $|V(T_{\hat{C_2}})|-1$ limbs of $T_{\hat{C_2}}[l,y]$ by non-decreasing height. Create a matrix $M:=M(T_{\hat{C_2}}[l,y], T_{C_2})$ having a row for each limb of $T_{\hat{C_2}}[l,y]$ (the rows ordered according to their corresponding height) and a column for each of the limbs of $T_{C_2}$. ($M$ plays the role of the limb embedding matrix.) 
		
		\item[II.2] For each limb $T_{\hat{C_2}}[a,b]$ of $T_{\hat{C_2}}[l,y]$ of height 1, and each limb $T_{C_2}[u,v]$ of $T_{C_2}$ set the entry associated with this pair of limbs in $M$ to 1 if either of the following two conditions is satisfied:
		\begin{itemize}
			\item[(i)] $a\in V(\hat{C_2})$, $u\in V(C_2)$ and $d_{T_{C_1}}(a)\geq b_{C_2}(u)$,
			\item[(ii)] $b\in V(\hat{C_2})$, $v\in V(C_2)$ and $d_{T_{C_1}}(b)\geq b_{C_2}(v)$.
		\end{itemize}
		Otherwise set the entry to 0.
		\item[II.3] For $h$ ranging between 2 and the height of $T_{\hat{C_2}}[l,y]$: \\For each limb $T_{\hat{C_2}}[a,b]$ of $T_{\hat{C_2}}[l,y]$ of height $h$, and each limb $T_{C_2}[u,v]$ of $T_{C_2}$, create a matrix $H:=H(T_{\hat{C_2}}[a,b],T_{C_2}[u,v])$, having a row for each neighbor of $b$ in $T_{\hat{C_2}}[a,b]-a$, and a column for each neighbor of $v$ in $T_{C_2}[u,v]-u$ ($H$ plays the role of the highest-limbs matrix.) For each neighbor $c$ of $b$,  $c\ne a$, and each neighbor $w$ of $v$, $w\ne u$, set the corresponding entry of $H$ to 1 if the $(T_{\hat{C_2}}[b,c],T_{C_2}[v,w])$-entry of $M$ is 1, and 0 otherwise. Determine a maximum bipartite matching for $H$. If it is complete for the rows and one of the conditions (i) or (ii) is satisfied, then set the value of the $(T_{\hat{C_2}}[a,b],T_{C_2}[u,v])$-entry of $M$ to 1. Otherwise set it to 0.

		\item[II.4] For each entry of $M$ which is 1 compute the isomorphism associated to it in the following way. If we have preserved the solution to all maximum matching problems for unit entries in $M$, then we can retrace, starting from the matching corresponding to this entry, the subtree of $T_{C_2}$ isomorphic to $T_{\hat{C_2}}[a,b]$ and the corresponding isomorphism. Note that it is not necessary to store these matchings as they can be computed again from $M$.
		
		If there is an entry in the last row of $M$ with value 1 such that the associated subtree of $T_{C_2}$, which is isomorphic to $T_{\hat{C_2}}$, contains all the cut-vertices of $C_2$, then return that there is a tree having $G$ as its exact-distance square.  Otherwise, return that there is no such a tree.
	\end{itemize}

	\noindent  To see that Stage I can be performed in polynomial time, it helps to bear in mind that the blocks and cut-vertices of a graph can be computed in quadratic time~\cite{10.1145/362619.362628}. Stage II can be performed in polynomial time since the number of bipartite matching problems we have to solve is at most quadratic in the number of vertices of $G$.
	
	Now we have to check the correctness of this algorithm. Before that we need a couple of lemmas.
	
	\begin{lemma}
		\label{restriction}
		Let $C$ and $D$ be two clique-trees. If $\varphi:T_{C}\rightarrow S\subseteq T_{D}$ is an isomorphism such that $\varphi(V(T_C)\cap V(C))\subseteq V(S)\cap V(D)$, then $\varphi(C)\subseteq V(D)$ and, moreover, the mapping $\psi: C\rightarrow D[\varphi(C)]$ defined by $\psi(v)=\varphi(v)$ is an isomorphism of the corresponding clique-trees.
	\end{lemma}
	
	\begin{proof}
		Let $v\in C$. By definition of $T_C$ we have that $v\in T_C$. By hypothesis we have that $\varphi(v)\in V(D)$.
		
		Now let us show that $\psi$ is an isomorphism. Clearly this function is a bijection. Let us see that it is an isomorphism, by showing that $uv\in E(C)$ is equivalent to $\varphi(u)\varphi(v)\in D$. It is clear that $uv\in E(C)$ if and only if there exists some $w\in V(T_C)\setminus V(C)$ such that $uwv$ is a path in $T_C$. But since $\varphi$ is an isomorphism this is also equivalent to having $\varphi(u)\varphi(w)\varphi(v)$ as a path in $S$ with $\varphi(w)\in S\setminus V(D)$. And in turn this is equivalent to having $\varphi(u)\varphi(v)\in D$, as desired.
	\end{proof}
	
	Let $G$ be a graph with two components $C_1$ and $C_2$ such that each of them is a clique-tree. Let $\hat{C_2}$ be the component of $\exact{T_{C_1}}{2}$ which is not $C_1$. We say a limb $T_{\hat{C_2}}[a,b]$ \emph{can be well embedded} in $T_{C_2}[u,v]$ if  $T_{\hat{C_2}}[a,b]$ can be embedded in $T_{C_2}[u,v]$ through an isomorphism which maps each vertex $s\in V(T_{\hat{C_2}}[a,b])\cap V(\hat{C_2})$ to a vertex $f_s\in V(C_2)$ such that $d_{T_{C_1}}(s) \ge b_{C_2}(f_s)$. We say that an isomorphism between $T_{\hat{C_2}}$ and a subtree of $T_{C_2}$ is a \emph{good isomorphism} if each vertex $s\in V(T_{\hat{C_2}})\cap V(\hat{C_2})$ is mapped to a vertex $f_s\in V(C_2)$ such that $d_{T_{C_1}}(s) \ge b_{C_2}(f_s)$. We say that a good isomophism is \emph{perfect} if the subtree of $T_{C_2}$ to which $T_{\hat{C_2}}$ is mapped contains all the cut-vertices of $C_2$.
	
	Our final ingredient is an analogue of Theorem~\ref{MatulaTheorem}.
	
	\begin{lemma}\label{ourmatching}
		Let $h\ge 2$ be an integer. Assume that the entries of $M$ have been computed in such a way that for a limb of $T_{\hat{C_2}}[l,y]$ of height $h-1$ and a limb of $T_{C_2}$, the position of $M$ corresponding to these limbs is $1$ if and only if  the limb of $T_{\hat{C_2}}[l,y]$ can be well embedded in that of $T_{C_2}$. A limb $T_{\hat{C_2}}[a,b]$ of $T_{\hat{C_2}}[l,y]$ of height $h$ can be well embedded to a limb $T_{C_2}[u,v]$ of $T_{C_2}$ if and only if the associated matrix $H:=H(T_{\hat{C_2}}[a,b],T_{C_2}[u,v])$ has a maximum bipartite matching which is complete for the rows and one of the conditions (i) and (ii) of II.2 is satisfied.
	\end{lemma}
	\begin{proof}
		First we assume that $T_{\hat{C_2}}[a,b]$ can be well embedded to $T_{C_2}[u,v]$. This implies that there is an isomorphism $\varphi$ from $T_{\hat{C_2}}[a,b]$ to $T_{C_2}[u,v]$ which is good and, in particular, maps $a$ to $u$. If we have $a\in V(\hat{C_2})$, then this tells us that $u\in V(C_2)$ and $d_{T_{C_1}}(a)\geq b_{C_2}(u)$, in other words, (i) is satisfied. If instead we have $b\in V(\hat{C_2})$, then we have that (ii) is satisfied.
		
		Recall that the entries of the matrix $H$ are filled using the information of the matrix $M$ for limbs of $T_{\hat{C_2}}[l,y]$ of height $h-1$. For each neighbor $c$ of $b$ in $T_{\hat{C_2}}[a,b]-a$ we have that $T_{\hat{C_2}}[b,c]$ can be embedded to $T_{C_2}[v,\varphi(c)]$. This and our assumption on $M$ implies that the $(c,\varphi(c))$-entry  of $H$ is 1. Since $\varphi$ is an isomorphism, we have that for each row in $H$ associated to a neighbor of $b$ in $T_{\hat{C_2}}[a,b]-a$ we have a different column in $H$ with a unit entry, i.e., $H$ has a maximum bipartite matching which is complete for the rows.

		For the other direction, notice that our assumption on $M$ implies that if $H$ has a maximum bipartite matching which is complete for the rows, then for each neighbor $c_i$ of $b$ in $T_{\hat{C_2}}[a,b]-a$ there exists a good isomorphism $\phi_i$ of $T_{\hat{C_2}}[b,c_i]$ to $T_{C_2}[v,\phi_i(c_i)]$. We define $\phi$ such that $\phi(a)=u$ and $\phi(x)=\phi_i(x)$ for every $x\in V(T_{\hat{C_2}}[b,c_i])$ and every neighbor $c_i$ of $b$ in $T_{\hat{C_2}}[a,b]-a$. Notice that $\phi$ is an isomorphism of $T_{\hat{C_2}}[a,b]$ to a limb of $T_{C_2}[u,v]$, in other words, $T_{\hat{C_2}}[a,b]$ can be embedded to $T_{C_2}[u,v]$. Moreover, if one of $(i)$ and $(ii)$ is satisfied, it follows that this isomorphism is good.
	\end{proof}

	\begin{theorem}
		Stages I and II correctly decide if the input graph $G$ admits an exact-distance square root which is a tree.
	\end{theorem}
	
	\begin{proof}
		We first show that, once Stage I has checked that $G$ has two components $C_1$ and $C_2$, each a clique-tree, it is enough for Stage II to decide if $T_{\hat{C_2}}$ has a perfect isomorphism to a subgraph of $T_{C_2}$. Indeed, if there is one such isomorphism $\varphi$, then by Lemma~\ref{restriction} (with $C=\hat{C_2}$ and $D=C_2$) we have an isomorphism $\psi\colon V(\hat{C_2}) \rightarrow C_{2}[\varphi(\hat{C}_2)]$ which, since $\varphi$ is perfect, satisfies the conditions of Theorem~\ref{fullcharac}, and guarantees the existence of a tree exact-distance square root. For the converse, if there is a tree exact-distance square root, then  Theorem~\ref{fullcharac} guarantees the existence of an isomorphism from $V(\hat{C_2})$ to $V(C_2)$ which can be easily extended to a perfect isomorphism from $T_{\hat{C_2}}$ to a subgraph of $T_{C_2}$.

		Now notice that deciding if $T_{\hat{C_2}}$ has a perfect isomorphism to a subgraph of $T_{C_2}$ is equivalent to deciding whether the limb $T_{\hat{C_2}}[l,y]$ rooted at a leaf $l$ of $T_{\hat{C_2}}$ can be well embedded in some limb of $T_{C_2}$ which contains all the cut-vertices of $C_2$.
		
		We will show, by induction on $h$, that steps $II.2$ and $II.3$ compute $M$ in such a way that for a limb $T_{\hat{C_2}}[a,b]$ of $T_{\hat{C_2}}[l,y]$ of height $h$ and a limb $T_{C_2}[u,v]$ of $T_{C_2}$, the position $(T_{\hat{C_2}}[a,b],T_{C_2}[u,v])$ of $M$ is 1 if and only if  $T_{\hat{C_2}}[a,b]$ can be well embedded in $T_{C_2}[u,v]$. Assuming this, and having computed all such entries, in step II.4 for each 1-entry of the last row the algorithm decides if all the cut-vertices are in the limb of $T_{C_2}[u,v]$. When that is the case, then we have an exact-distance square root of $G$ which is a tree. Otherwise there is no such a root. This verifies the correctness of the algorithm.

		To prove the claim by induction, first for the base of induction, let $T_{\hat{C_2}}[a,b]$ be a limb of $T_{\hat{C_2}}[l,y]$ of height $h=1$. Note that $T_{\hat{C_2}}[a,b]$ can be embedded in any limb of $T_{C_2}$, and that one of the conditions $II.2.(i)$ and $II.2.(ii)$ is satisfied if and only if $T_{\hat{C_2}}[a,b]$ can be well embedded in the limb of $T_{C_2}$. 
		
		For the main step of induction, assume the entries of $M$ have been computed preserving the desired property for all rows corresponding to limbs of $T_{\hat{C_2}}[l,y]$ of height at most $h-1$ for some $h\geq2$.
		Let $T_{\hat{C_2}}[a,b]$ be a limb of $T_{\hat{C_2}}[l,y]$ of height $h$ and let $T_{C_2}[u,v]$ be a limb of $T_{C_2}$. Consider the matrix $H:=H(T_{\hat{C_2}}[a,b],T_{C_2}[u,v])$ associated to these limbs. Note that the limbs $T_{\hat{C_2}}[b,c]$ of $T_{\hat{C_2}}[a,b]$, where $c\ne a$ is a neighbor of $b$, are also limbs of $T_{\hat{C_2}}[l,y]$ having height at most $h-1$. Similarly, the limbs $T_{C_2}[v,w]$ of $T_{C_2}$, where $w\ne u$ is a neighbor of $v$, are also limbs of $T_{C_2}$. Hence, the information needed to fill the entries of $H$ is already available in the previously computed portion of $M$. Now step  $II.3$ determines a maximum matching for $H$ and, by Lemma \ref{ourmatching}, $T_{\hat{C_2}}[a,b]$ can be well embedded in $T_{C_2}[u,v]$ if and only if the matching is complete for the rows and one of conditions (i) or (ii) (of II.2) is met. Since this step puts a 1 in the entry of $M$ corresponding to $T_{\hat{C_2}}[a,b]$ and $T_{C_2}[u,v]$ if and only if  the matching is complete for the rows  and one of conditions (i) and (ii) are met, this entry of $M$ is computed in a way that preserves the desired property. \end{proof}

	At the start of this proof we saw that once our algorithm has decided that there is a perfect isomorphism from $T_{\hat{C_2}}$ to a subgraph of $T_{C_2}$, then from this isomorphism we can obtain, through Lemma~\ref{restriction}, an isomorphism satisfying the conditions of Theorem~\ref{fullcharac}. From this isomorphism, we can obtain a corresponding tree exact-distance square root by using the procedure given in the proof of Theorem \ref{nearcharac}.

	\section{Tree exact-distance square roots are not unique}\label{sec:unique}

	In this section, we prove Theorem~\ref{theorem:unique}, i.e., for every $m\ge 2$ we give a graph with at least $m!$ nonisomorphic exact-distance square roots each of which is a tree. 
	
	\begin{figure}[ht]
		\centering
		\includegraphics[height=2.2 in]{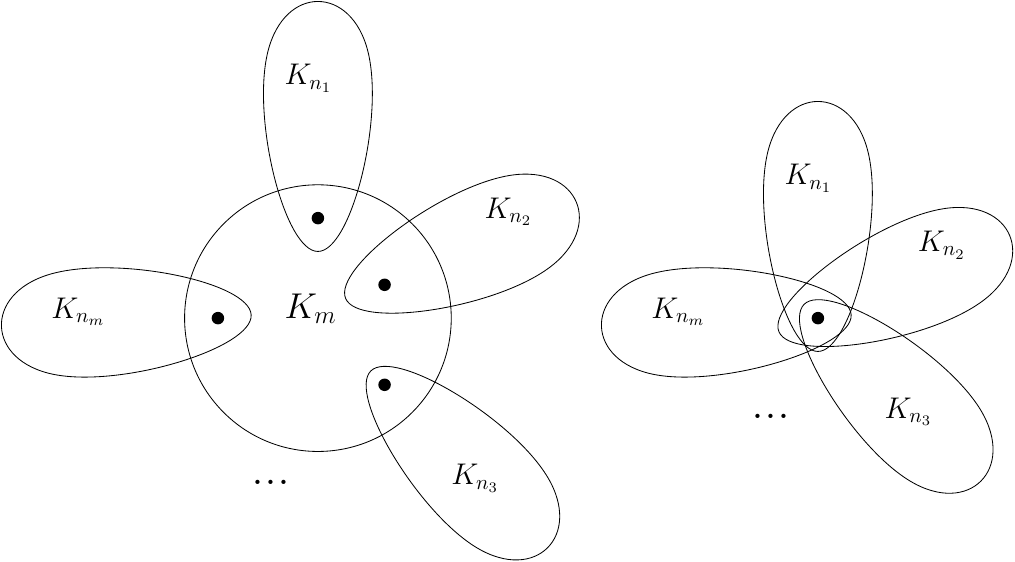}
		\caption{The graph $G_{S}$}
		\label{twoflowers}
	\end{figure}
	
	\begin{figure}[ht]
		\centering
		\includegraphics[height=1.5 in]{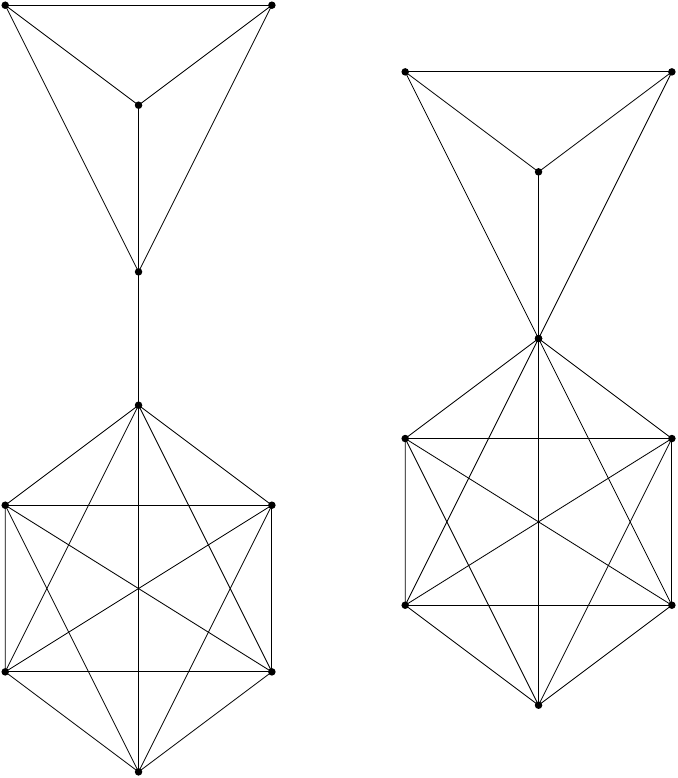}
		\caption{$G_{S}$ for the sequence $S=(4,6)$}
		\label{S46flowers}
	\end{figure}

	\begin{proof}[Proof of Theorem~\ref{theorem:unique}]
		Let $S=(n_1,n_2,\ldots,n_m)$ be a sequence of $m$ positive integers such that $1<n_1<n_2<\ldots<n_m$. We define the graph $G_{S}$ with two connected components, as follows. For the first component we take a copy of $K_m$, with vertices $v_1,\ldots,v_m$, and a copy of each of $K_{n_1-1},\ldots,K_{n_m-1}$. Then, for each $1\le i\le m$, we make $v_i$ adjacent to all the vertices of $K_{n_i-1}$ (forming  a $K_{n_i}$). The other component is formed by taking a copy of the first component and contracting all the vertices $v_1,\ldots,v_m$, into one. See Figures \ref{twoflowers} and \ref{S46flowers} for illustrations.
		
		\begin{figure}[h!]
			\centering
			\includegraphics[height=2.2 in]{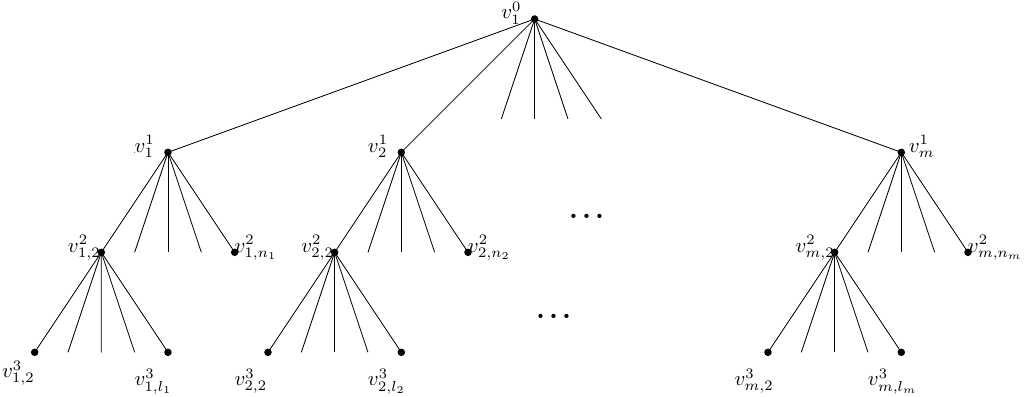}
			\caption{The tree $T_L$}
			\label{rootofflowers}
		\end{figure}

		\begin{figure}[h!]
			\centering
			\includegraphics[angle=90, height=1 in]{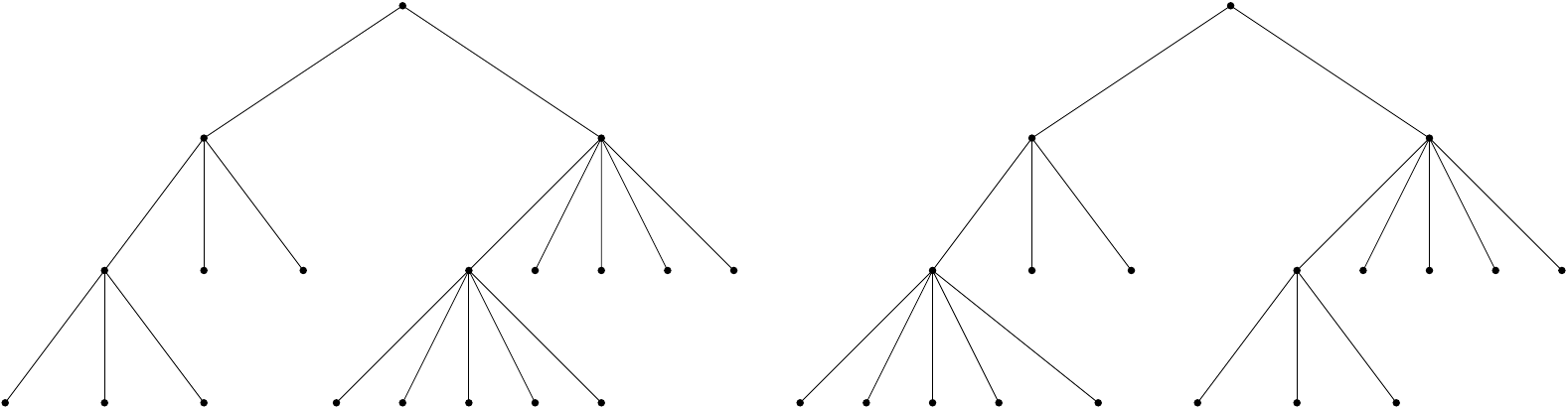}
			\caption{Two trees which have $G_{S}$, the graph from Figure~\ref{S46flowers}, as their exact-distance square}
			\label{S46rootof}
		\end{figure}
		
		Now we claim that there exist at least $m!$ nonisomorphic	trees all of which have $G_S$ as  their exact-distance square. Let $L=(l_1,\ldots,l_m)$ be a permutation of $n_1,n_2,\ldots,n_m$; we are going to construct a tree for this permutation. First, take a vertex $v_{1}^{0}$ with $m$ neighbours $v_{1}^{1},\ldots,v_{m}^{1}$. Then for each $i\in\{1,\ldots,m\}$, we make $v_{i}^{1}$ adjacent to $n_i-1$ new vertices $v_{i,2}^{2},\ldots,v_{i,n_1}^{2}$. Finally, to $v_{i,1}^{2}$ we connect $l_i-1$ new vertices $v_{i,2}^{3},\ldots,v_{i,l_i}^{3}$. Let $T_L$ be the tree obtained (see Figure~\ref{rootofflowers}). 
		It is not hard to see that the exact-distance square of $T_L$ is isomorphic to $G_S$, and that if $L'\neq L$ is a permutation of $n_1,n_2,\ldots,n_m$, then the graphs $T_L$ and $T_{L'}$ are not isomorphic.
	\end{proof}

	\medskip
	{\bf Acknowledgment.} We would like to thank Pierre Le Bodic for discussions on an early stage of this project. The project has received funding from the following grants. FONDECYT/ANID Iniciaci\'on en Investigaci\'on Grant 11201251, Programa Regional MATH-AMSUD MATH210008,   ANR-France project HOSIGRA (ANR-17-CE40-0022), National Natural Science Foundation of China (Grant Nos. 12131013, 12242111, 11601430), Guangdong Basic \& Applied Basic Research Foundation (Grant Nos. 2023A1515030208, 2022A1515010899), and Shaanxi Fundamental Science Research Project for Mathematics and Physics (Grant No. 22JSZ009).

	\bibliographystyle{plain}
	\bibliography{ref}

\end{document}